\documentclass{amsart}

\usepackage{amsmath}
\usepackage{amscd}
\usepackage{amssymb}

\input xy
\xyoption{all}

\newcommand{\bga}{\boldsymbol{\gamma}} \newcommand{\bpsi}{\bar{\psi}}
\newcommand{\pd}{\partial}
\newcommand{\bC}{{\mathbb C}} \newcommand{\bc}{{\bf c}}
\newcommand{\bE}{{\mathbb E}} \newcommand{\bF}{{\mathbb F}}

\newcommand{\bP}{{\mathbb P}}

\newcommand{\bZ}{{\mathbb Z}}

\newcommand{\cB}{{\mathcal B}}

\newcommand{\cH}{{\mathcal H}}

\newcommand{\cM}{{\mathcal M}}
\newcommand{\cO}{{\mathcal O}}

\newcommand{\half}{\frac{1}{2}}

\newcommand{\cX}{{\mathcal X}}

\newcommand{\Mbar}{\overline{\cM}}

\newcommand{\cl}[1]{[\![{#1} ]\!]}
\newcommand{\cor}[1]{\langle {#1} \rangle}

\DeclareMathOperator{\Aut}{Aut} \DeclareMathOperator{\ad}{ad}
\DeclareMathOperator{\diag}{diag} \DeclareMathOperator{\End}{End}
\DeclareMathOperator{\ev}{ev}
\DeclareMathOperator{\Hom}{Hom} 
 \DeclareMathOperator{\rank}{rank}
\DeclareMathOperator{\class}{Class} 
\DeclareMathOperator{\ch}{ch}

\newtheorem{theorem}{Theorem}[section]
\newtheorem{theorem/definition}{Theorem/Definition}[section]

\newtheorem{proposition}{Proposition}[section]
\newtheorem{lemma}{Lemma}[section]

\theoremstyle{remark}
\newtheorem{remark}{Remark}[section]

\theoremstyle{definition}
 \newtheorem{example}{Example}[section]
\newtheorem{definition}{Definition}[section]

\newcommand{\be}{\begin{equation}}
\newcommand{\ee}{\end{equation}}
\newcommand{\bea}{\begin{eqnarray}}
\newcommand{\ben}{\begin{eqnarray*}}
\newcommand{\een}{\end{eqnarray*}}
\newcommand{\eea}{\end{eqnarray}}
\newcommand{\bet}{\begin{equation}
\begin{split}}
\newcommand{\eet}{\end{split}
\end{equation}}

\begin{document}

\title{On Computations of Hurwitz-Hodge Integrals}
\author{Jian Zhou}
\address{Department of Mathematical Sciences\\
Tsinghua University\\Beijing, 100084, China}
\email{jzhou@math.tsinghua.edu.cn}

\begin{abstract}
We describe a method to compute Hurwitz-Hodge integrals.
\end{abstract}
\maketitle

\section{Introduction}

Recently there have been a lot of interests in the Gromov-Witten theory of
the orbifold $[\bC^3/\bZ_3]$,
both from physicists and mathematicians.
In physics, $[\bC^3/\bZ^3]$ represents the orbifold point in the A-model moduli space for the local $\bP^2$.
By mirror symmetry,
Aganagic, Bouchard and Klemm \cite{Aga-Bou-Kle} studied its mirror $B$-model and made
some predictions on the Gromov-Witten invariants for $[\bC^3/\bZ_3]$.
In mathematics,
it is interesting to verify these predictions and extend the results to other orbifolds.

The mathematical theory of orbifold Gromov-Witten invariants has been
developed for symplectic orbifolds by Chen-Ruan \cite{Che-Rua}
and for Deligne-Mumford stacks by Abramovich-Graber-Vistoli \cite{Abr-Gra-Vis}.
Bryan and Graber \cite{Bry-Gra} introduced the notion of twisted degrees so that
the potential function for orbifolds can be suitably defined.
The physicists' predictions has been verified in various case by different authors:
Coates-Corti-Iritani-Tseng \cite{Coa-Cor-Iri-Tse},
Bayer-Cadman \cite{Bay-Cad}, Cadman-Cavalieri \cite{Cad-Cav},
Bouchard-Cavalieri \cite{Bou-Cav}.
The problem in general can be phrased as the computations of Hurwitz-Hodge integrals
on the moduli spaces $\Mbar_{g, n}(\cB G)$ of twisted stable maps to the classifying
stack of a finite group $G$.
In this paper we describe a method that can be used to  compute
Hurwitz-Hodge integrals and hence can be used to compute Gromov-Witten invariants of
other orbifolds.
As examples,
we present some details for the orbifold $[\bC^3/\bZ_5(1,1,3)]$.
In principle it is easy to automate this method and we are working on a Maple program
for $G=\bZ_3$.

Our approach is in the same spirit as the recursive calculations of ordinary Hodge integrals.
Recall ordinary $\psi$-integrals on the Deligne-Mumford moduli spaces can
be computed using the famous Witten-Kontsevich theorem \cite{Wit, Kon}.
Faber \cite{Fab} described an algorithm that computes  Hodge integrals recursively by reducing to the $\psi$-integrals
and he implemented it by a beautiful Maple program.
This is based on Mumford's GRR relations \cite{Mum},
which have been generalized to moduli spaces of stable maps by Faber-Pandhariande \cite{Fab-Pan}
and they derived a system of differential equations for the generating function
of Hodge integrals.
Givental \cite{Giv1} gave a solution to this system in terms of the generating function
of $\psi$-integrals,
and he reformulated the results in an enlightening quantization formulation in a later work \cite{Giv2}.
This formulation was generalized by Coates-Givental \cite{Coa-Giv} to Hodge integral type
invariants (called twisted invariants) on moduli spaces of stable maps to projective manifolds,
and further generalizations to orbifolds was made by Tseng \cite{Tse}.
As has already been allured to in Tseng's paper,
our starting point is  that by specializing his work to the classifying stack $\cB G$,
the computations of Hurwitz-Hodge integrals are reduced
to the generating functions of $\psi$-integrals on $\Mbar_{g, n}(\cB G)$ of a finite group $G$,
and the latter has been computed by Jarvis-Kimura \cite{Jar-Kim}.
It is then possible to extend Faber's algorithm to compute Hurwitz-Hodge integrals.
We carry out some calculations by hand and leave the implementation of  this algorithm for future work.
Of course most of the materials in this paper are already in the literature,
but surprisingly they have not been put together to compute Hurwitz-Hodge integrals.
In the $G=\bZ_2$ case we recover a formula due to Faber-Padharipande \cite{Fab-Pan2};
in the $G=\bZ_3$ case we verify some of the predictions of Aganagic-Bouchard-Klemm \cite{Aga-Bou-Kle};
and in the $\bZ_5$ case we present some examples which seem to be new.

\section{Orbifold Cohomology of $\cB G$ and Quantization of Its Formal Loop Space}

In this section we recall some natural bases of the orbifold cohomology of the classifying stack $\cB G$ of a finite group $G$.
We also recall Givental quantization of its formal loop space following Tseng \cite{Tse}.

\subsection{Some bases of the orbifold cohomology of $\cB G$}

Geometrically,
the classifying stack $\cB G$ of a finite group is a point with a trivial $G$-action.
Its inertia orbifold is $I \cB G = \coprod_{[\![\gamma]\!]} \cB G_{[\![\gamma]\!]}$,
where $\cl{\gamma}$ denotes the conjugacy class of $\gamma \in G$,
and $G_{\cl{\gamma}}$ denotes its centralizer in $G$.

The orbifold cohomology of $\cB G$ is, as a vector space,
$$\cH := H^*_{orb}(\cB G,\bC) :=
H^*(I\cB G,\bC) = \bigoplus_{[\![\gamma]\!]} \bC.$$
For each conjugacy class $[\![\gamma]\!]$ in $G$,
let $e_{[\![\gamma]\!]}$ denote the class $1 \in H^0(\cB G_{[\![\gamma]\!]}, \bC) \subseteq \cH$.
The $e_{[\![\gamma]\!]}$'s form a basis of $\cH$ which we will call the {\em class basis}.
As a special case of the ring structure on orbifold cohomology introduced by Chen-Ruan \cite{Che-Rua},
there is a structure of a Frobenius algebra on $\cH$,
where the multiplication is given by:
\be
e_{\cl{\gamma_1}}e_{\cl{\gamma_2}} = \sum_{\sigma_i \in \cl{\gamma_i}, i=1, 2}
\frac{|C(\sigma_1\sigma_2)|}{|G|} e_{\cl{\sigma_1\sigma_2}},
\ee
and the metric is given by:
\be
\eta(e_{\cl{\gamma_1}}, e_{\cl{\gamma_2}}) = \frac{1}{|C(\gamma_1)|} \delta_{\cl{\gamma_1}\cl{\gamma_2^{-1}}}.
\ee
See Jarvis-Kimura \cite{Jar-Kim} for details.
For later use,
let $\eta_{\cl{\gamma_1},\cl{\gamma_2}} = \eta(e_{\cl{\gamma_1}}, e_{\cl{\gamma_2}})$
and let $(\eta^{\cl{\gamma_1},\cl{\gamma_2}})$ be the inverse matrix of
$(\eta_{\cl{\gamma_1},\cl{\gamma_2}})$.

The Frobenius algebra $\cH$ is isomorphic to the Frobenius algebra $\class_{\bC}(G)$ of class functions of $G$,
where $e_{[\![\gamma]\!]}$ corresponds to  the characteristic function $\chi_{\cl{\gamma}}$
of the class $\cl{\gamma}$.
On $\class_{\bC}(G)$ the product is the convolution product $\star$ of functions:
\ben
\chi_{\cl{\gamma_1}} \ast \chi_{\cl{\gamma_2}}
= \sum_{\substack{\sigma_1, \sigma_2 \\ \sigma_i \in \cl{\gamma_i}}}
\frac{|C(\sigma_1 \sigma_2)|}{|G|} \chi_{\cl{\sigma_1 \sigma_2}},
\een
and the metric is:
$$\langle f, g \rangle
= \frac{1}{|G|} \sum_{\sigma \in G} f(\sigma)g(\sigma^{-1}).$$
Denote by $Z\bC(G)$ the center of the group algebra of $G$.
On $\bC[G]$ we have the the standard group-algebra product $\cdot$,
and a metric
$$(\alpha,\beta) =
\frac{1}{|G|}\delta_{\alpha,\beta^{-1}}.$$
When restricted to $Z\bC(G)$ they provide a structure of a Frobenius algebra.
Let $\Phi$ denote the standard additive isomorphism $\Phi: \class_{\bC}(G) \to Z\bC[G]$
defined by linearly extending
$\Phi(\chi_{\cl{\gamma}}) = \sum_{\alpha \in \cl{\gamma}} \alpha$.
Then $\Phi$ is an isomorphism of Frobenius algebras.

With these isomorphisms,
one can form another natural basis of $\cH$.
Let $\{ V_\alpha\}_{\alpha=1}^r$ be the set of irreducible
representations of $G$ and let $\chi_\alpha$ denote the character
of $V_\alpha$.  For all $\alpha=1,\ldots,r$, the elements
\[
f_\alpha := \frac{\dim V_\alpha}{|G|}\sum_{g\in G}
\chi_{\alpha}(g^{-1}) g
\]
form a basis of $Z\bC[G]$ and satisfy the following equations:
\bea
&& f_{\alpha_1} \cdot f_{\alpha_2} =
\delta_{\alpha_1,\alpha_2} f_{\alpha_1}, \\
&& (f_{\alpha_1},f_{\alpha_2}) =
\delta_{{\alpha_1},{\alpha_2}} \nu_{\alpha_1},
\eea
where for all $\alpha=1, \ldots, r$,
\[
\nu_\alpha = \left(\frac{\dim V_\alpha}{|G|}\right)^2.
\]
Furthermore,
the identity element satisfies:
$$1 = \sum_{\alpha=1}^r f_{\alpha}.$$
The corresponding basis in $\cH$,
also denoted by $\{f_{\alpha}\}$,
will be called the {\em representation basis}.
It is easy to see that
\be \label{eqn:EtoF}
f_{\alpha} = \frac{\dim V_\alpha}{|G|}\sum_{\cl{\gamma}}
\chi_{\alpha}(\cl{\gamma^{-1}}) e_{\cl{\gamma}}
\ee
Using the orthogonality relation
\ben
\sum_{\alpha} \chi_{\alpha}(g^{-1})\chi_{\alpha}(h)
= \begin{cases}
|G_g|, & \cl{g} = \cl{h}, \\
0, & \text{otherwise},
\end{cases}
\een
one gets:
\be \label{eqn:FtoE}
e_{\cl{\gamma}} = \sum_{\alpha} \frac{|G/G_{\gamma}|}{\dim V_{\alpha}} \chi_{\alpha}(\cl{\gamma}) f_{\alpha}.
\ee

\subsection{Givental quantization on the formal loop space of a Frobnius algebra}

Let $H$ be a $\bC$-Frobenius algebra of dimension  $h$ with metric
$\langle \cdot, \cdot \rangle$.
One can define a canonical symplectic form $\Omega$ on $H[z,z^{-1}]$ by:
\be
\Omega(f(z), g(z)) = \frac{1}{2\pi i} \oint \langle f(-z), g(z)\rangle dz.
\ee
Fix a pair of dual bases $\{v_1, \dots, v_h\}$, $\{w_1, \dots, w_h\}$ of $(H, \langle\cdot, \cdot \rangle$.
I.e.,
$$\langle v_i, w_j \rangle = \delta_{ij},$$
for $i, j =1, \dots, h$.
Let $g_{ij}= \langle v_i, v_j\rangle$,
and let $g^{ij}$ denote the inverse matrix.
Then one has $w_j = \sum_i v_i g^{ij}$.
For $f(z) \in H[z, z^{-1}]$,
write
\ben
f(z) = \sum_{k=0}^{\infty} \sum_{i=1}^h p_k^i w_i (-z)^{-k-1}
+ \sum_{k=0}^{\infty} \sum_{i=1}^h q_k^i v_i z^k.
\een
Then $\{p^i_k, q^i_k\}_{i=1, \dots, h, k \geq 0}$
are Darboux coordinates on $(H[z, z^{-1}], \Omega)$.
The associated Fock space is the space of polynomials in $\{ q^i_k\}_{i=1, \dots, h, k \geq 0}$
and an extra formal variable $\lambda$.
On this space define operators
\begin{align*}
\widehat{q^i_k} : & = \lambda^{-1} \cdot q^i_k \cdot, &
\widehat{p^i_k}: & = \lambda \frac{\pd}{\pd q^i_k}.
\end{align*}
The operators $\{\widehat{q^i_k}\}$ are called the {\em creators},
and the operator $\{\widehat{p^i_k}\}$ are called the {\em annihilators}.
Given a polynomial $f$ in $\{p^i_k, q^i_k\}$,
the quantization of $f$ is the operator $\hat{f}$ obtained from $f$
by replacing $q^i_k$ with $\widehat{q^i_k}$,
$p^i_k$ by $\widehat{p^i_k}$,
and apply the normal ordering, i.e., put all the annihilators on the right of the creators.
E.g.,
$$:\widehat{p^i_k}\widehat{q^j_l}: = :\widehat{q^j_l}\widehat{p^i_k}:
=\widehat{q^j_l}\widehat{p^i_k}.$$

An {\em infinitesimal symplectic transformation} is an element $A \in \End H[z, z^{-1}]$
such that
\be
\Omega(Af, g) = \Omega(f, Ag).
\ee
Given an infinitesimal symplectic transformation $A$,
one can associate a quadratic polynomial  $P_A$ by:
$$P_A(f) = \half \Omega(Af, f).$$
We are interested in $\hat{A}:=\widehat{P_A}$.

\begin{example} \cite{Coa} \label{ex:Coates}
Let $B \in \End(H)$, $Bv_\alpha = \sum_\beta B_{\; \alpha}^\beta v_\beta$.
As usual,
set $B_{\alpha\beta} = \sum_\gamma g_{\alpha\gamma}B^\gamma_{\;\beta}$
and $B^{\alpha\beta} = \sum_{\gamma} B^\alpha_{\;\,\gamma} g^{\gamma\beta}$.
Suppose that $A=Bz^m$ is an infinitesimal symplectomorphism.
That is to say,
if $m$ is odd, $B$ is self-adjoint;
if $m$ is even, $B$ is anti-self-adjoint.
Set
$$\pd_{\alpha,k} = \frac{\pd}{\pd q^\alpha_k}.$$
A straightforward calculation shows that
if $m < 0$ then
\be
\hat{A} = \lambda^{-2} \sum_k (-1)^{k+m} B_{\alpha\beta}q^\beta_kq^\alpha_{-1-k-m}
- \sum_k B^{\alpha}_{\;\,\beta} q^\beta_k\pd_{\alpha,k+m},
\ee
and if $m \geq 0$ then
\be
\hat{A} = - \sum_{k} B^\alpha_{\;\,\beta} q_k^\beta\pd_{\alpha, k+m}
+ \frac{\lambda^2}{2} \sum_{k } (-1)^k B^{\alpha\beta}\pd_{\beta,k}\pd_{\alpha,m-1-k}.
\ee
\end{example}

\subsection{Infinitesimal symplectic transformations on $\cH$
induced by vector bundles on $\cB G$}

A representation $V_{\rho}$ of $G$ determines an orbifold vector bundle $E_{\rho}$ on $\cB G$,
and also an orbifold vector bundle $I E_{\rho}$ on the inertia orbifold $I\cB G$.
Denote by $E_{\rho, \cl{\gamma}}$ the restriction of $I E_{\rho}$ on $\cB G_{\cl{\gamma}}$.
Denote by $r:=|\gamma|$ the order of $\gamma$ in $G$.
Then we have a decomposition
$$V_{\rho} = \bigoplus_{0 \leq l \leq r} V_{\rho, \gamma}^{(l)},$$
where $\gamma$ acts on $V_{\rho,\gamma}^{(l)}$ with eigenvalue $\zeta_r^l$ and $\zeta_r=\exp(2\pi i/r)$
is a primitive $r$-th root of unity.
Denote by $E_{\rho, \cl{\gamma}}^{(l)}$ the subbundle of $E_{\rho, \cl{\gamma}}$
corresponding to $V_{\rho, \gamma}^{(l)}$.
Denote by $E_{\rho, \cdot}^{(l)}$ the bundle on $I\cB G$
whose restriction to the sector $\cB G_{\cl{\gamma}}$ is $E^{(l)}_{\rho, \cl{\gamma}}$.

For a multiplicative characteristic class $\bc(\cdot) :=\exp (\sum_{k=0}^{\infty} s_k \ch_k(\cdot))$
of vector bundles,
the $(\bc, E_\rho)$-twisted inner product on $\cH$ defined in \cite{Tse}
is given by:
\ben
\langle e_{\cl{\gamma_1}}, e_{\cl{\gamma_2}} \rangle_{\bc(F)}
= \frac{1}{|C(\gamma_1)|} \exp (s_0 \dim V_{\rho,\gamma_1}^{(0)}) \delta_{\cl{\gamma_1}\cl{\gamma_2^{-1}}}.
\een
By \cite[Lemma 3.2.1]{Tse},
the map $(\cH, \langle \cdot, \cdot \rangle_{\bc(E_{\rho})}) \to (\cH, \langle\cdot, \cdot\rangle)$
defined by:
$a \mapsto a\sqrt{\bc(E_{\rho, \cdot}^{(0)})}$ is an isometry.
Here the multiplication is the ordinary multiplication given by {\em componentwise} multiplications
on $H^*(\cB G_{\cl{\gamma}})$.
The determinant of this map is denoted by $\det \sqrt{\bc(E_{\rho, \cdot}^{(0)})}$.
Because
\ben
\sqrt{\bc(E_{\rho, \cdot}^{(0)})}
& = & \sum_{\cl{\gamma}} \exp (\half s_0 \dim E_{\rho, \gamma}^{(0)} ) e_{\cl{\gamma}},
\een
therefore the ordinary multiplication by $\sqrt{\bc(E_{\rho, \cdot}^{(0)})}$
is given by the diagonal matrix $\diag (\exp (\half s_0 \dim E_{\rho, \gamma}^{(0)}))$
in the $\{e_{\gamma}\}$ basis,
hence we have
\be
\det\sqrt{\bc(E_{\rho, \cdot}^{(0)})}
= \prod_{\cl{\gamma}} \exp (\half s_0 \dim E_{\rho, \gamma}^{(0)} ).
\ee

Recall that the Bernoulli polynomials $B_m(x)$ are defined by
$$\frac{te^{tx}}{e^t-1}=\sum_{m\geq 0}\frac{B_m(x)t^m}{m!}.$$
The Bernoulli numbers $B_m$ are given by $B_m:=B_m(0)$.

\begin{definition}
For each integer $m\geq 0$, define an element $A_m(E_{\rho}) \in H^*(I\cB G)
=\oplus_{\cl{\gamma}}H^*(\cB G_{\cl{\gamma}})$ as follows:
$$A_m(E_{\rho})|_{\cB G_{\cl{\gamma}}}:=\sum_{0\leq l\leq |\gamma|-1}
(\dim V_{\rho, \gamma}^{(l)}) B_m(l/|\gamma|) e_{\cl{\gamma}}.$$
\end{definition}

As a special case of \cite[Corollary 4.1.5]{Tse},
we have

\begin{lemma}
Ordinary multiplications by the following classes define infinitesimally symplectic transformations
on $(\cH,\Omega)$ and $(\cH_{(\bc,E_{\rho})},\Omega_{\bc(E_{\rho})})$:
$$A_{2m}(E_{\rho})z^{2m-1},\,\, A_{2m+1}(E_{\rho})z^{2m},\,\,\,m\geq 1;\,\, A_0(E_{\rho})/z,\,\,
A_1(E_{\rho})+\frac{1}{2}\ch_0(q^*F^{(0)}).$$
\end{lemma}

\section{The stack of orbifold stable maps into $\cB G$}

We recall some basic facts about $\Mbar_{g, n}(\cB G)$.

\subsection{The stack $\Mbar_{g, n}(\cB G)$}

Let $\Mbar_{g, n}(\cB G)$ be the stack of $n$-pointed orbifold stable maps into $\cB G$,
in the sense of Chen and Ruan \cite{Che-Rua}.
They are called balanced twisted stable maps by Abramovich and Vistoli \cite{Abr-Vis}.
The stack $\Mbar_{g,n}(\cB G)$ is a smooth, proper Deligne-Mumford stack of dimension $3g- 3 + n$
with projective coarse moduli space \cite[Thm 3.0.2]{Abr-Cor-Vis}.

To understand these spaces better,
let us first recall the description of an orbifold stable map from a smooth, $n$-pointed
orbicurve into $\cB G$  \cite{Jar-Kim}.
The domain is a smooth $n$-pointed orbicurve $(\Sigma, x_1,\dots, x_n)$
which has non-trivial orbifold structure only at marked points $x_1, \dots, x_n$.
An orbifold stable map from such an orbicurve  into $\cB G$ is a smooth curve $\widetilde{\Sigma}$
with a $G$-action and a $G$-equivariant map
$p: \widetilde{\Sigma} \to \Sigma$,
such that $p': \widetilde{\Sigma}' \to \Sigma'$ is principal $G$ bundle,
where $\Sigma' = \Sigma-\{x_1, \dots, x_n\}$,
$\widetilde{\Sigma}' = \widetilde{\Sigma} -  p^{-1}(\{x_1, \dots, x_n\})$,
and $p':=p|_{\widetilde{\Sigma}}$.
The stabilizer $G_{x_i}$ of a marked point $x_i$ of $\Sigma$ is always cyclic,
and the holonomy around $x_i$ determines an injective map $G_{x_i} \to G$.
The inverse image of $x_i$ is an orbit of $G$ isomorphic to $G/G_{x_i}$,
i.e.,
$p^{-1}(x_i)$ consists of $|G/G_{x_i}|$ points $y_{i,1}, \dots, y_{i,|G/G_{x_i}|}$.
It follows that the Euler numbers of $\widetilde{\Sigma}$ and $\Sigma$ are related by:
\be \label{eqn:Euler}
\chi(\widetilde{\Sigma}) - \sum_{i=1}^n |G/G_{x_i}|
= |G| (\chi(\Sigma) - n).
\ee

A general stable orbifold map $f:\Sigma \to \cB G$ is also
a $G$-equivariant map $f: \widetilde{\Sigma} \to \Sigma$,
but now both $\widetilde{\Sigma}$ and $\Sigma$ are allowed to be nodal curves.
We say a $G$-orbit of the $G$-action on $\widetilde{\Sigma}$ is {\em free}
if it consists of $|G|$ points.
The nonfree orbits are only allowed to be the inverse images under $p$ of
 the marked points and the nodes   of $\Sigma$.
If $y \in \Sigma$ is a node such that $p^{-1}(y)$ is nonfree,
then for any $q \in p^{-1}(y)$,
$q$ is a node of $\widetilde{\Sigma}$,
and the stabilizer $G_q$ is a cyclic group $\bZ_l$.
Furthermore,
let $C_1$ and $C_2$ be two branches that meet at $q$,
then the $G_q$-representation on $T_qC_1$ is the conjugate representation of the $G_q$-representation
on $T_qC_2$.
Let $\Sigma'$ and $\widetilde{\Sigma}'$ denote the regular parts of $\Sigma$
and $\widetilde{\Sigma}$, respective.
Then $p':=p|_{\widetilde{\Sigma}'}: \widetilde{\Sigma}' \to \Sigma'$ is a principal $G$-bundle.

A homomorphism of two such maps $p_i: \tilde{\Sigma}_i \to \Sigma_i$ ($i=1, 2$)
is a commutative diagram of holomorphic maps:
\ben
\xymatrix{
  \widetilde{\Sigma}_1 \ar[d]_{p_1} \ar[r]^{\tilde{f}}
                & \widetilde{\Sigma}_2 \ar[d]^{p_2}  \\
  \Sigma_1 \ar[r]_{f}
                & \Sigma_2            }
\een
This also defines the automorphism group $\Aut(p: \widetilde{\Sigma} \to \Sigma)$.
We have an exact sequence
\ben
1 \to G \to \Aut(p: \widetilde{\Sigma} \to \Sigma) \to \Aut(\Sigma) \to 1.
\een
If the automorphism group is finite,
then we say the orbifold holomorphic map $f: (\Sigma, x_1, \dots, x_n) \to \cB G$ is stable.

\subsection{Structure forgetting morphism}

By the above discussions we have seen that
if $p: (\widetilde{\Sigma}, \{y_{i,j}\}_{1 \leq i \leq n, 1 \leq j \leq |G/G_{x_i}|}) \to (\Sigma, x_1, \dots, x_n)$
is stable if and only if $(\Sigma, x_1, \dots, x_n)$ is a stable curve,
and so we have a {\em forgetting morphism}
$\phi: \Mbar_{g, n}(\cB G) \to \Mbar_{g, n}$ by ``forgetting" the orbifold structure on $\Sigma$.
We will call this morphism the structure forgetting morphism.

The principal $G$-bundle $p: \widetilde{\Sigma}' \to \Sigma'$
determines a class of homomorphisms $\pi_1(\Sigma') \to G$ up to conjugations by $G$.
This class is invariant under automorphisms of $p: \widetilde{\Sigma} \to \Sigma$.
Conversely,
 a class of such homomorphisms  determine a unique class of principal $G$-bundles on $\Sigma'$.

Because $\pi_1(\Sigma')$ is generated by
$\{a_i, b_i\}_{i=1, \dots, g} \cup \{ c_j \}_{j=1, \dots, n}$,
subject to the relation
\be \label{eqn:Monodromy}
\prod_{i=1}^g [a_i, b_i] = \prod_{j=1}^n c_j.
\ee
Fix conjugacy classes $\bga_{[n]} = (\cl{\gamma_1}, \dots, \cl{\gamma_n})$,
let
$\Hom (\pi_1(\Sigma'), G)(\bga_{[n]})$
be the homomorphisms whose images of $c_j$ lies in the conjugacy class $\gamma_j$,
$j=1, \dots, n$.
Define
\ben
\cX^G_{g} (\bga_{[n]})& := & \{(\alpha_1, \dots, \alpha_g, \beta_1, \dots,
\beta_g, \sigma_1, \dots, \sigma_n) \in G^{2g+n} |\textstyle \\
&& \prod^g_{i=1}[\alpha_i,\beta_i] =\prod^n_{j=1}\sigma_j, \ \sigma_j
\in [\![\gamma_j]\!]\text{ for all } j \},
\een
and let $G$ acts diagonally on $\cX^G_{g}(\bga_{[n]})$ by conjugations.
Then we have an 1-1 correspondence
$$\Hom (\pi_1(\Sigma'), G)(\bga_{[n]})/\ad G \cong \cX^G_{g}(\bga_{[n]})/\ad G.$$
Thus we have the following:

\begin{proposition}
For a given curve $[C, p_1,\dots,p_n]$ corresponding to a point of the smooth locus $\cM_{g,n}$,
the fiber of $\Mbar_{g,n}(\cB G)$ over the point $[C]$ corresponds to the quotient
$\operatorname{Hom}(\pi_1(C-\{p_1,\dots,p_n\}, G)/ ad G$.
\end{proposition}

It follows that$\cM_{g, n}(\cB G, , \cl{\gamma_1},\dots,\dots,\cl{\gamma_n}) \to \cM_{g, n}$
is a finite morphism of degree
$$\Omega^G_{g}(\bga) =\frac{|\cX^G_{g} (\bga)|}{|G|}.$$
These numbers can be computed as follows

\begin{lemma} \cite{Jar-Kim}
The numbers $\Omega^G_g({\bf\gamma})$ depend only on the conjugacy
classes $[\![\gamma_i]\!]$, are independent of the ordering of the
$\gamma_i$ in $\bga$, and satisfy the following relations:
\begin{enumerate}
\item  \textbf{Cutting trees:}  For $g=g_1+g_2$ and $I \coprod J=
\{1,\dots,n\}$, let $\bga_I = (\gamma_{i_1}, \dots,
\gamma_{i_{|I|}})$ and $\bga_J=(\gamma_{j_1}, \dots,
\gamma_{j_{|J|}})$ $$\Omega^G_g (\bga) =
\Omega^G_{g_1}(\bga_I,\zeta) \eta^{[\![\zeta]\!][\![\xi]\!]}
\Omega^G_{g_2}(\xi,\bga_J)$$
\item \textbf{Cutting loops:} $$\Omega^G_g(\bga) =
\eta^{[\![\zeta]\!](\xi]\!]} \Omega^G_{g-1}(\zeta,\xi,\bga)$$
\item \textbf{Forgetting tails:} $$\Omega^G_g (\bga) = \Omega^G_g
(1,\bga)$$
\end{enumerate}
\end{lemma}

When $G$ is abelian, the computation of $\Omega^G_g(\cl{\gamma}_{[n]})$ is much simpler.
Note we have
\be
\cX^G_g(\cl{\gamma}_{[n]}) = \begin{cases}
|G|^{2g}, & \gamma_1 \cdots \gamma_n  =1, \\
\emptyset, & \gamma_1 \cdots \gamma_n  \neq 1.
\end{cases}
\ee
Hence we get
\be
\Omega^G_g(\cl{\gamma}_{[n]}) = \begin{cases}
|G|^{2g-1}, & \gamma_1 \cdots \gamma_n  =1, \\
0, & \gamma_1 \cdots \gamma_n  \neq 1.
\end{cases}
\ee

\subsection{Other natural morphisms}

When the holonomy $\gamma_{n+1}$  around the marked point $x_{n+1}$ is trivial,
we may \emph{forget} the data of that marked point.
This gives a morphism \cite{Abr-Vis}:
$$\pi: \Mbar_{g,n+1}(\cB G, [\![\gamma]\!]_{[n]},[\![1]\!]) \to
\Mbar_{g, n}(\cB G, [\![\gamma]\!]_{[n]}).$$
The marked points define natural sections
$$\sigma_1, \dots, \sigma_n: \Mbar_{g, n}(\cB G, [\![\gamma]\!]_{[n]})
\to  \Mbar_{g,n+1}(\cB G, [\![\gamma]\!]_{[n]},[\![1]\!]).$$
Denote by  $D_{i, \cl{\gamma}_{[n]}}$ the corresponding divisors,
$i=1, \dots, n$.

For each $I=1, \dots, n$,
there is an evaluation map $\ev_i: \Mbar_{g, n}(\cB G, [\![\gamma]\!]_{[n]}) \to
\cB G_{\cl{\gamma_i}}$
For a conjugacy class $[\![\gamma]\!]$ of $G$,
one has the following two types of degree $2$ morphisms:
$$\Mbar_{g,n+2}(\cB G, [\![\gamma]\!]_{[n]},[\![\gamma]\!], [\![\gamma^{-1}]\!]) \to
\Mbar_{g, n}(\cB G, [\![\gamma]\!]_{[n]}),$$
$$\Mbar_{g_1,n_1+1}(\cB G, [\![\gamma_I]\!],[\![\gamma]\!]) \times
\Mbar_{g_2,n_2+1}(\cB G, [\![\gamma_J]\!],[\![\gamma^{-1}]\!]) \to
\Mbar_{g, n}(\cB G, [\![\gamma]\!]_{[n]}),$$
where $I \cup J = [n]$.

\section{Intersection Theory of $\Mbar_{g, n}(\cB G)$}

In this section we recall some results of Jarvis-Kimura \cite{Jar-Kim},
and some results of Tseng \cite{Tse} specialized to $\cB G$.

\subsection{Intersection theory of psi classes on $\Mbar_{g, n}(\cB G)$}

Following \cite{Jar-Kim} and \cite{Tse},
define
$$\bpsi_i: = \phi^*\psi_i$$
where $\phi: \Mbar_{g, n}(\cB G) \to \Mbar_{g,n}$ is the structure forgetting map.
For $\beta_1, \dots, \beta_n \in \cH$,
define the \emph{$n$-point correlators} by
$$\langle \tau_{a_1}(\beta_1)\dots \tau_{a_n}(\beta_n)\rangle^G_g
:= \int_{\Mbar_{g,n}(\cB G)} \prod_i \bpsi_i^{a_i} \ev_i^*(\beta_i).$$
Let ${\bf t} = \{ t^{\cl\gamma}_a\}$ be formal variables associated to the class basis $\{e_{\cl\gamma}\}$.
The \emph{potential function} of the correlators is defined by
$F^G({\bf t}) = \sum_g F_g({{\bf t}}) \lambda^{2g-2}$ in $\lambda^{-2}\bC[[\bf t,\lambda]]$,
where $F^G_g(\bf t) := \cor{\exp(\bf t\cdot\bf \tau)}^G_g$,
where ${\bf t}\cdot{\bf \tau} = \sum_{a=0}^{\infty} \sum_{\cl\gamma} t_a^{\cl\gamma} \tau_a(e_{\cl\gamma})$.
Let $Z^G := \exp(F^G)$.
When $G$ is trivial,
we write $Z^G$ and $F^G$ as $Z$ and $F$ respectively.
They are invariant of a point.

Let $\bf u$ be formal variables $\{u_a^\alpha\}$ for all integers $a\geq 0$
and $\alpha=1,\ldots,r$ associated to the representation basis $\{f_\alpha\}$.
For each $\alpha = 1,\ldots,r$, let
$\tilde{\bf u}^\alpha$ be formal variables $\{\,\tilde{u}_a^\alpha\,\}$,
where $a\geq 0$ and $\tilde{u}_a^\alpha := (\nu_\alpha)^{\frac{1-a}{3}} u_a^\alpha$.
By (\ref{eqn:EtoF}) and (\ref{eqn:FtoE}), we have
\be \label{eqn:TtoU}
u^\alpha = \sum_{\cl\gamma} \frac{|G/G_{\gamma}|}{\dim V_{\alpha}}
\chi_{\alpha}(\cl{\gamma}) t^{\cl\gamma},
\ee
\be \label{eqn:UtoT}
t^{\cl\gamma} = \sum_{\alpha} \frac{\dim V_\alpha}{|G|}
\chi_{\alpha}(\cl{\gamma^{-1}}) u^\alpha.
\ee

\begin{proposition} \cite{Jar-Kim}
We have
\begin{equation}
F^G({\bf t}) = \sum_{\alpha=1}^r F(\tilde{\bf u}^\alpha).
\end{equation}
\end{proposition}

From this it is straightforward to deduce the Virasoro constraints
and KdV hierarchy equations for $F^G$.
For details,
see \cite{Jar-Kim}.
This Proposition follows from the following:

\begin{proposition} \cite{Jar-Kim}
The correlators in the class basis are related to the usual correlators by:
\be
\langle\tau_{a_1}(\cl{\gamma_1})\dots\tau_{a_n}(\cl{\gamma_n})\rangle^G_g
= \langle\tau_{a_i}\dots \tau_{a_n}\rangle_g \Omega^G_g (\bga).
\ee
The correlators in the representation basis are related to the usual coorelators by:
\begin{equation}
\cor{\tau_{a_1}(f_{\alpha_1})\cdots\tau_{a_n}(f_{\alpha_n})}^G_g =
\begin{cases}
\nu_\alpha^{1-g} \cor{\tau_{a_1}\cdots\tau_{a_n}}_g, & \alpha_1 = \ldots = \alpha_n = \alpha, \\
0, & \text{otherwise}.
\end{cases}
\end{equation}
\end{proposition}

\subsection{The $J$-function of $\cB G$}

Recall  the orbifold $J$-function $J_{\cB G}(t,z)$ of $\cB G$ is defined by \cite{Tse}:
$$J_{\cB G}(t,z)=z+t+\sum_{n \geq 3} \frac{1}{(n-1)!}\sum_{k\geq 0,\,\alpha}
\cor{\tau_0(t)^{n-1}\tau_k(\phi_\alpha)}^G_0\frac{\phi^\alpha}{z^{k+1}}.$$
Here $\{\phi_\alpha\}$ is an additive basis of $\cH$ and $\{\phi^\alpha\}$ is its dual basis
under the orbifold pairing $(\cdot,\cdot)_{orb}$.
Let $\{\phi_\alpha\}$ be the representation basis $\{f_{\alpha}\}$,
and let $t = \sum_\alpha u^\alpha f_\alpha$,
then we have:
\ben
J_{\cB G}(t,z) & = & z + \sum_\alpha u^\alpha f_\alpha
+ \sum_{n \geq 3} \frac{1}{(n-1)!}\sum_{k\geq 0,\,\alpha}
\cor{\tau_0(\sum_\alpha u^\alpha f_\alpha)^{n-1}\tau_k(f_\alpha)}^G_0\frac{f_\alpha}{\nu_\alpha z^{k+1}} \\
& = & z+ z \sum_\alpha f_\alpha \biggl(u^\alpha z^{-1}
+ \sum_{n \geq 3} \frac{(u^\alpha)^{n-1}}{(n-1)!} z^{-(n-1)} \biggr) \\
& = & z + z \sum_\alpha f_{\alpha} (e^{u^\alpha/z}-1)
\een

\subsection{Hurwitz-Hodge integrals}

Because $R^0\pi _{*}f^{*} (E_{\alpha })$ is the trivial vector bundle corresponding to $V_{\alpha}^G$,
one gets a vector bundle $R^1\pi _{*}f^{*} (E_{\alpha})$.
Define
$$\bF^\alpha_{g,n} = R^0\pi _{*}f^{*} (E_{\alpha }) - R^1\pi _{*}f^{*} (E_{\alpha })
\in K(\Mbar_{g,n}(\cB G)).$$
Its virtual rank is defined by:
$$\rank \bF^\alpha_{g,n} = \rank R^0\pi _{*}f^{*} (E_{\alpha }) - \rank R^1\pi _{*}f^{*} (E_{\alpha }).$$

\begin{proposition}
On $\Mbar_{g, n}(\cB G; \cl{\gamma}_{[n]})$ the virtual rank of $\bF_{g, n}^{\alpha}$ is given by:
\be \label{eqn:Rank}
\rank \bF^\alpha_{g,n}  = (\dim V_\alpha)(1-g) - \sum_{i=1}^n  \sum_{l=0}^{|\gamma_i|}
(\dim V_{\alpha, \gamma_i}^{(l)})\frac{l}{|\gamma_i|}.
\ee
\end{proposition}

\begin{proof}
On $\Mbar_{g, n}(\cB G; \cl{\gamma}_{[n]})$,
by the orbifold Riemann-Roch formula of Kawasaki \cite{Kaw} we have:
\ben
&& \rank R^0\pi _{*}f^{*} (E_{\alpha }) - \rank R^1\pi _{*}f^{*} (E_{\alpha }) \\
& = & \half (\dim V_{\alpha}) \int_{\Sigma} c_1(T\Sigma) \\
& + & \sum_{i=1}^n \frac{1}{|\gamma_i|} \sum_{l=0}^{|\gamma_i|} \dim V_{\alpha, \gamma}^{(l)} \sum_{j=1}^{|\gamma_i|-1}
\frac{e^{2\pi \sqrt{-1} jl/|\gamma_i|}}{1 - e^{-2\pi j \sqrt{-1} /|\gamma_i|}},
\een
where on the right-hand side,
the second term is the contribution from the twisted sector
and can be computed using Lemma \ref{lm:Sum} below,
and the first term is the contribution from the main sector, where $T\Sigma$ stands for
orbifold tangent bundle of $\Sigma$.
Therefore,
\ben
&& \int_{\Sigma} c_1(T\Sigma)
= \frac{1}{|G|} \int_{\widetilde{\Sigma}} c_1(T\widetilde{\Sigma})
= \frac{1}{|G|} \chi (\widetilde{\Sigma})
= \chi(\Sigma) - \sum_{i=1}^n (1 - \frac{|G/\langle\gamma_i\rangle|}{|G|})\\
& = & \chi(\Sigma) - \sum_{i=1}^n (1 - \frac{1}{|\gamma_i|}).
\een
Here we have used (\ref{eqn:Euler}).
Hence
\ben
\rank \bF^\alpha_{g,n}  & = & \frac{\dim V_\alpha}{2}(2-2 g - \sum_{i=1}^n (1 - \frac{1}{|\gamma_i|}))
 + \sum_{i=1}^n  \sum_{l=0}^{|\gamma_i|} \frac{\dim V_{\alpha, \gamma_i}^{(l)}}{|\gamma_i|}
(\frac{|\gamma_i|-1}{2} - l) \\
& = & (\dim V_\alpha)(1-g) - \sum_{i=1}^n  \sum_{l=0}^{|\gamma_i|}
(\dim V_{\alpha, \gamma_i}^{(l)})\frac{l}{|\gamma_i|}.
\een
\end{proof}

\begin{lemma} \label{lm:Sum}
For $m > 1$ and $0 \leq l < m$, we have
\be
\sum_{j=1}^{m-1} \frac{\xi_m^{jl}}{1 - \xi_m^{-j}}
= \frac{m-1}{2} - l,
\ee
where $\xi_m = e^{2\pi i /m}$.
\end{lemma}

\begin{proof}
First recall
\ben
\sum_{j=1}^{m-1} \xi_m^{aj} =
m \delta_{a,0}- 1.
\een
Notice that
\ben
&& \sum_{j=1}^{m-1} \xi_m^{jl} \frac{z^{m-1} + z^{m-2} + \cdots + 1}{z- \xi_m^{-j}} \\
& = & \sum_{j=1}^{m-1} \xi_m^{jl} (z^{m-2} + (1 + \xi_m^{-j}) z^{m-3} +
(1 + \xi_m^{-j} + \xi_m^{-2j}) z^{m-4} + \cdots \\
&& +(1 + \xi_m^{-j} + \cdots + \xi_m^{-(m-2)j})).
\een
Now take $z=1$ on both sides:
\ben
&& m \sum_{j=1}^{m-1} \frac{\xi_m^{jl}}{1 - \xi_m^{-j}}
=  \sum_{j=1}^{m-1} \xi_m^{jl} \sum_{k=0}^{m-2} (m-1-k)\xi_m^{-jk} \\
& = & \sum_{k=0}^{m-2} (m-1-k) (m \delta_{k, l} - 1)
= \frac{m(m-1)}{2} - ml.
\een
This finishes the proof.
\end{proof}

Define the $\alpha$-twisted Hurwitz-Hodge classes by:
\be
\lambda_{j, \alpha} = (-1)^j c_j(R^1\pi^*f^*E_{\alpha}).
\ee
We also define
\be
\ch_{k, \alpha} = \ch_k(R^1\pi^*f^*E_{\alpha}).
\ee
By the standard relationship between Chern classes and Chern characters,
each $\lambda_{j, \alpha}$ can be expressed as a polynomial in $\ch_{k, \alpha}$ and vice versa.

By a {\em Hurwitz-Hodge integral} we mean an integral of the form
\be
\int_{\Mbar_{g, n}(\cB G, \cl{\gamma})} \prod_{i=1}^n \psi_i^{a_i} \cdot
\prod_{j=1}^m \lambda_{b_j, \alpha_j},
\ee
where $\alpha_1, \dots, \alpha_m$ are $G$-characters.

Fix  a multiplicative characteristic class $\bc (\cdot)
= \exp (\sum_{k=0}^{\infty} s_k \ch_k(\cdot))$ and a $G$-character $\alpha$.
Following \cite{Tse},
define {\em $(\bc,E_\alpha)$-twisted total descendant potential} of $\cB G$ by:
$$Z_{(\bc,E_{\alpha})}({\bf t}):
= \exp{\left(\sum_{g=0}^{\infty} \lambda^{2g-2} \sum_{n \geq 1}
\frac{1}{n!}\int_{\Mbar_{g,n}(\cB G)} \bc(\bF^\alpha_{g,n})
\prod_{i=1}^n \sum_{k=1}^{\infty}
\ev_i^*(\sum_{\cl\gamma} t^{\cl\gamma}_ke_{\cl\gamma}) \bpsi_i^k \right)}.$$
Now we state Tseng's formula \cite{Tse} applied to our special case:

\begin{theorem}
We have
\be \label{eqn:Tseng}
Z_{(\bc,E_{\alpha})}({\bf t})
= (\det\sqrt{\bc(E_{\alpha})})^{\frac{1}{24}}
\exp\left(\sum_{k \geq 0} s_k \frac{A_{k+1}(E_{\alpha})z^{k}}{(k+1)!}
+s_0\frac{\dim E_{\alpha, \cdot}^{(0)}}{2}\right)^\wedge Z^G.
\ee
\end{theorem}

In this formula,
the variables $q^{\cl\gamma}_k$ are related to $t^{\cl\gamma}_k$ by the dilaton shift:
$$q^{\cl\gamma}_k = t^{\cl\gamma}_k - \delta_{k1}\delta_{\cl{\gamma} \cl1}.$$

\begin{remark}
In the above formula,
we have used only one $G$-character.
There is a similar formula when one has more than one $G$-character.
This is because the formula is proved by converting a Mumford type GRR relation to a system of differential equations.
\end{remark}

Because $Z^G$ has already been calculated by Jarvis-Kimura \cite{Jar-Kim},
this Theorem provides an effective way to compute Hurwitz-Hodge integrals.

\begin{example}
Take derivative in $s_k$ on both sides of (\ref{eqn:Tseng})
and then set all $s_i = 0$:
\ben
&& \sum_{g=0}^{\infty} \lambda^{2g-2} \sum_{n \geq 1}
\frac{1}{n!}\int_{\Mbar_{g,n}(\cB G)} \ch_{k,\alpha} \prod_{i=1}^n \sum_{l=1}^{\infty}
\ev_i^*(\sum_{\cl\gamma} t^{\cl\gamma}_l e_{\cl\gamma}) \bpsi_i^l \cdot Z^G \\
& = & \left(\frac{A_{k+1}(E_{\alpha})z^{k}}{(k+1)!}\right)^\wedge Z^G
\een
By Example \ref{ex:Coates},
the operator
$\cO_k(\alpha): = \left(\frac{A_{k+1}(E_{\alpha})z^{k}}{(k+1)!}\right)^\wedge$ is of the form:
$$\cO_k(\alpha)
= D_1 + \frac{\lambda^2}{2} \sum_{l=0 }^{k-1} C^{\alpha\beta}\pd_{\beta,l}\pd_{\alpha,k-1-l},$$
where $D_1$ is a first order differential operator,
and $\pd_{\beta,l} = \frac{\pd}{\pd t^{\cl\beta}_l}$.
Therefore,
\ben
&& (Z^G)^{-1} \left(\frac{A_{k+1}(E_{\alpha})z^{k}}{(k+1)!}\right)^\wedge Z^G \\
& = & (Z^G)^{-1} (D_1 + \frac{\lambda^2}{2} \sum_{l=0 }^{k-1} C^{\alpha\beta}\pd_{\beta,l}\pd_{\alpha,k-1-l})
\exp F^G \\
& = & D_1 F^G + \frac{\lambda^2}{2} \sum_{l=0 }^{k-1} C^{\alpha\beta}\pd_{\beta,l}\pd_{\alpha,k-1-l}F^G
+ \frac{1}{2} \sum_{l=0 }^{k-1} C^{\alpha\beta}\pd_{\beta,l} F^G \cdot\pd_{\alpha,k-1-l}F^G.
\een
Recall $F^G= \sum_{g=0}^{\infty} \lambda^{2g-2} F_g^G$.
Hence by taking the coefficients of $\lambda^{2g-2}$,
we get:
\be \label{eqn:F0}
\begin{split}
&\sum_{n \geq 3}
\frac{1}{n!}\int_{\Mbar_{0,n}(\cB G)} \ch_{k, \alpha} \prod_{i=1}^n
\sum_{l=1}^{\infty} \ev_i^*(\sum_{\cl\gamma} t^{\cl\gamma}_le_{\cl\gamma}) \bpsi_i^l  \\
= & D_1 F^G_0
+ \frac{1}{2} \sum_{l=0 }^{k-1} C^{\alpha\beta}\pd_{\beta,l} F^G_0
\cdot\pd_{\alpha,k-1-l}F^G_0,
\end{split}
\ee
and for $g > 0$,
\be \label{eqn:Fg}
\begin{split}
&\sum_{n > 2-2g}
\frac{1}{n!}\int_{\Mbar_{g,n}(\cB G)} \ch_{k, \alpha}
\prod_{i=1}^n \sum_{l=1}^{\infty} \ev_i^*(\sum_{\cl\gamma} t^{\cl\gamma}_le_{\cl\gamma}) \bpsi_i^l  \\
 = & D_1 F^G_g
+ \frac{1}{2} \sum_{l=0 }^{k-1} C^{\alpha\beta}\pd_{\beta,l}\pd_{\alpha,k-1-l}F^G_{g-1} \\
 + & \frac{1}{2} \sum_{h=0}^g \sum_{l=0 }^{k-1} C^{\alpha\beta}\pd_{\beta,l} F^G_h
\cdot\pd_{\alpha,k-1-l}F^G_{g-h}.
\end{split}
\ee
\end{example}

In the next three sections we will present some examples in the case of $G=\bZ_2$, $\bZ_3$ and $\bZ_5$.

\section{Examples: The Case of $\cB Z_2$}

\subsection{Orbifold cohomology of $\cB Z_2$}

Let $\omega$ be a generator of $\bZ_2$.
In the class basis,
the product on the orbifold cohomology of $\cB \bZ_2$ is given by:
\be
e_{\cl{\omega^a}}e_{\cl{\omega^b}} =
 e_{\cl{\omega^{a+b}}},
\ee
where $a, b \in \bZ_2$,
and the metric is given by:
\be
\eta(e_{\cl{\omega^a}}, e_{\cl{\omega^b}}) = \begin{cases}
\frac{1}{2}, & a+ b \cong 0 \pmod{2}, \\
0, & \text{otherwise}.\end{cases}
\ee
I.e.,
\be
(\eta_{\alpha\beta})
= \begin{pmatrix}
\frac{1}{2} & 0 \\ 0 & \frac{1}{2}
\end{pmatrix}
\ee
therefore,
\be
(\eta^{\alpha\beta})
= \begin{pmatrix}
2 & 0 \\ 0 & 2
\end{pmatrix}
\ee

\subsection{Some Hurwitz-Hodge integrals in the $g=0$ case}
Denote by $V_0$ the trivial $1$-dimensional representation of $\bZ_2$,
and by $V_1$ the nontrivial $1$-dimensional representation of $\bZ_2$.
We are interested in the vector bundle $E_1$ on $\cB \bZ_2$ associated with $V_1$.
By the monodromy condition (\ref{eqn:Monodromy}),
$\Mbar_{g, n}(\cB \bZ_2; \cl1^k, \cl\omega^{n-k})$ is nonempty only if $n-k$ is an even number
$2m \geq 0$.
By the rank formula (\ref{eqn:Rank}),
\be
\rank R^1\pi_*f^*V_1 = m+g-1
\ee
on $\Mbar_{g, k+2m}(\cB \bZ_2; \cl1^k, \cl\omega^{2m})$.

Now we can recover a result due to Faber-Pandharipande \cite{Fab-Pan2}.

\begin{proposition}
For $m > 1$ we have
\ben
\frac{1}{(2m)!}\int_{\Mbar_{0,2m}(\cB \bZ_2; \cl\omega^{2m})} \ch_{2m-3}(\bF^1_{0,2m})
=\frac{B_{2m-2}}{(2m-2)!} (2^{2m-2}-1).
\een
\end{proposition}

\begin{proof}
It is straightforward to see that
\ben
&& \left(\frac{A_{k+1}(E_1)z^{k}}{(k+1)!} \right)^\wedge  \\
& = & - \sum_{l=0}^{\infty} \big(\frac{B_{k+1}}{(k+1)!} q_l^0\pd_{0, k+l}
+ \frac{B_{k+1}(1/2)}{(k+1)!} q_l^1\pd_{1, k+l}\big) \\
& + & \frac{1}{2} \hbar \big(2\frac{B_{k+1}(0)}{(k+1)!} \sum_{l=0}^{k-1} (-1)^l \pd_{0,l}\pd_{0,k-1-l}
+ 2 \frac{B_{k+1}(1/2)}{(k+1)!} \sum_{l=0}^{k-1} (-1)^l \pd_{1,l} \pd_{1,k-1-l}\big).
\een
Therefore,
\ben
&& \sum_{n \geq 3}
\frac{1}{n!}\int_{\Mbar_{0,n}(\cB \bZ_2)} \ch_{k, 1}  \prod_{i=1}^n
\sum_{k=1}^{\infty} \ev_i^*(t^{0}_k e_{\cl1} + t^1_k e_{\cl\omega}) \bpsi_i^k  \\
& =  &  - \sum_{l=0}^{\infty} \big(\frac{B_{k+1}}{(k+1)!} q_l^0\pd_{0, k+l}F_0
+ \frac{B_k(1/2)}{k!} q_l^1\pd_{1, k+l} F_0^{\bZ_2} \big) \\
& + & \frac{1}{2} \big(2\frac{B_{k+1}(0)}{(k+1)!} \sum_{l=0}^{k-1} (-1)^l (\pd_{0,l} F_0^{\bZ_2} )(\pd_{0,k-1-l}F_0) \\
& + & 2 \frac{B_{k+1}(1/2)}{(k+1)!} \sum_{l=0}^{k-1} (-1)^l (\pd_{1,l}F_0^{\bZ_2}) (\pd_{1,k-1-l}F_0^{\bZ_2})\big).
\een
Now take $k = 2m-3$ and consider the coefficients of $(t^1_0)^{2m}$.
On the left-hand side we get
\ben
\frac{1}{(2m)!}\int_{\Mbar_{0,2m}(\cB \bZ_2; \cl\omega^{2m})} \ch_{2m-3,1}
\een
On the right-hand side we get:
\ben
&& \frac{B_{2m-2}}{(2m-2)!} \cor{\tau_0(e_{\cl\omega})^{2m} \tau_{2m-2}(e_{\cl1})}_0^{\bZ_2} \\
& - & \frac{B_{2m-2}(1/2)}{(2m-2)!} \cor{\tau_0(e_{\cl\omega})^{2m-1}\tau_{0}(e_{\cl1})\tau_{2m-2}(e_{\cl1})}_0^{\bZ_2} \cdot 2m \\
& + & \frac{B_{2m-2}}{(2m-2)!} \sum_{0 \leq l \leq 2m-4; l even}
(2m)! \frac{1}{(l+2)!} \cor{\tau_0(e_{\cl\omega})^{l+2}\tau_l(e_{\cl1})}_0^{\bZ_2} \\
&& \cdot \frac{1}{(2m-2-l)!} \cor{\tau_0(e_{\cl\omega})^{2m-2-l}\tau_{0}(e_{\cl1})\tau_{2m-4-l}(e_{\cl1})}_0^{\bZ_2} \\
& - & \frac{B_{2m-2}(1/2)}{(2m-2)!} \sum_{0 \leq l \leq 2m-4; l odd}
(2m)! \frac{1}{(l+2)!} \cor{\tau_0(e_{\cl\omega})^{l+2}\tau_l(e_{\cl\omega})}_0^{\bZ_2} \\
&& \cdot \frac{1}{(2m-2-l)!} \cor{\tau_0(e_{\cl\omega})^{2m-2-l}\tau_{2m-4-l}(e_{\cl\omega})}_0^{\bZ_2} \\
& = & \frac{B_{2m-2}}{(2m-2)!} \frac{1}{2}
- \frac{B_{2m-2}(1/2)}{(2m-2)!} \frac{1}{2} \cdot 2m \\
& + & \frac{B_{2m-2}}{(2m-2)!} \sum_{0 \leq l \leq 2m-4; l even}
\binom{2m}{l+2}  \frac{1}{2^2}
- \frac{B_{2m-2}(1/2)}{(2m-2)!} \sum_{0 \leq l \leq 2m-4; l odd}
\binom{2m}{l+2} \frac{1}{2^2} \\
& = & \frac{1}{4} \frac{B_{2m-2}}{(2m-2)!}
\sum_{0 \leq l \leq 2m; l even}
\binom{2m}{l}
- \frac{1}{4} \frac{B_{2m-2}(1/2)}{(2m-2)!} \sum_{0 \leq l \leq 2m; l odd}
\binom{2m}{l} \\
& = & \frac{1}{4} \frac{B_{2m-2}}{(2m-2)!} 2^{2m-1}
- \frac{1}{4} \frac{B_{2m-2}(1/2)}{(2m-2)!} 2^{2m-1} \\
& = & \frac{B_{2m-2}}{(2m-2)!} (2^{2m-2}-1).
\een
In the last equality we have used the identity:
\ben
B_n(1/2) = (\frac{1}{2^{n-1}}-1)B_n,
\een
which can be proved as follows.
\ben
\sum_{n \geq 0} \frac{B_n(1/2)}{n!} t^n = \frac{te^{t/2}}{e^t-1}
= 2 \frac{t/2}{e^{t/2} -1} - \frac{t}{e^t-1}
= \sum_{n \geq 0} (\frac{1}{2^{n-1}}-1)\frac{B_n}{n!} t^n.
\een
This completes the proof.
\end{proof}

\section{Example: The Case of $\cB \bZ_3$}

\subsection{Orbifold cohomology of $\cB Z_3$}

Let $\omega$ be a generator of $\bZ_3$.
In the class basis,
the product on the orbifold cohomology of $\cB \bZ_3$ is given by:
\be
e_{\cl{\omega^a}}e_{\cl{\omega^b}} =
 e_{\cl{\omega^{a+b}}},
\ee
and the metric is given by:
\be
\eta(e_{\cl{\omega^a}}, e_{\cl{\omega^b}}) = \begin{cases}
\frac{1}{3}, & a+ b \cong 0 \pmod{3}, \\
0, & \text{otherwise}.\end{cases}
\ee
I.e.,
\be
(\eta_{\alpha\beta})
= \begin{pmatrix}
\frac{1}{3} & 0 & 0 \\ 0 & 0 & \frac{1}{3} \\ 0 & \frac{1}{3} & 0
\end{pmatrix}
\ee
and so
\be
(\eta^{\alpha\beta})
= \begin{pmatrix}
3 & 0 & 0 \\ 0 & 0 & 3 \\ 0 & 3 & 0
\end{pmatrix}
\ee

\subsection{Potential function of $\cB \bZ_3$}

Let $\{ V_\alpha\}_{\alpha=0}^2$ be the set of irreducible
representations of $\bZ_3$,
where $\omega$ acts on $V_{\alpha}$ as multiplication by $\xi_3^\alpha$.
The representation basis is given by:
\[
f_\alpha := \frac{1}{3}\sum_{a=0,1,2}
\xi_3^{-a\alpha}e_{\cl{\omega^a}}.
\]
One has
\be
e_{\cl{\omega^a}} = \sum_{\alpha=0}^2 \xi_3^{a\alpha} f_{\alpha}.
\ee
Let $\{t^0, t^1, t^2\}$ be linear coordinates in $\{e_{\cl1}, e_{\cl\omega}, e_{\cl{\omega^2}}\}$
and let $\{u^0, u^1, u^2\}$ be linear coordinates in $\{f_0, f_1, f_2\}$.
We have
\begin{align}
t^a & = \frac{1}{3}\sum_{\alpha=0}^2 \xi_3^{-a\alpha} u^\alpha, &
u^{\alpha} & = \sum_{a=0}^2 \xi_3^{a\alpha} t^a.
\end{align}
Hence
\ben
F^{Z_3}({\bf u}) = \sum_{\alpha=0}^2 F(\tilde{{\bf u}}^{\alpha}),
\een
where
$\tilde{\bf u}_a^\alpha = 3^{2(a-1)/3} u_a^\alpha$.
This formula is very effective to automatically compute $F^{\bZ_3}$.
For example,
using Maple we get:
\ben
F^{\bZ_3}_0 & = & \frac{1}{3}t^0_0t^1_0t^2_0+\frac{1}{18}(t^0_0)^3+\frac{1}{18}(t^1_0)^3
+\frac{1}{18}(t^2_0)^3 \\
& + & \frac{1}{6} (t^1_0)^2t^2_0t^2_1+\frac{1}{18}(t^0_0)^3t^0_1+\frac{1}{3}t^0_0t^1_0t^2_0t^0_1
+\frac{1}{6}(t^0_0)^2t^1_0t^2_1+\frac{1}{6}(t^0_0)^2t^2_0t^1_1 \\
&+& \frac{1}{6}t^0_0(t^1_0)^2t^1_1
+\frac{1}{6}t^0_0(t^2_0)^2t^2_1+\frac{1}{6}t^1_0(t^2_0)^2t^1_1+\frac{1}{18}(t^1_0)^3t^0_1
+\frac{1}{18}(t^2_0)^3t^0_1 + \cdots, \\
F^{\bZ_3}_1 & = & \frac{1}{8}t^0_1 + \frac{1}{8}t^2_0t^1_2+\frac{1}{8}t^1_0t^2_2
+\frac{1}{8}t^0_0t^0_2+\frac{1}{8}t^1_1t^2_1+\frac{1}{16}(t^0_1)^2 + \cdots
\een

\subsection{Some Hodge integrals}
Because of the monodromy condition (\ref{eqn:Monodromy}),
the component $\Mbar_{g,n}(\cB \bZ_3; \cl1^{n_0}, \cl\omega^{n_1}, \cl{\omega^2}^{n_2})$ is nonempty only if
$n_1 \equiv n_2 \pmod{3}$.
Let $\bE^\vee_a = R^1\pi_*f^*V_a$.
By the rank formula (\ref{eqn:Rank}),
we have on the component $\Mbar_{g,n}(\cB \bZ_3; \cl1^{n_0}, \cl\omega^{n_1}, \cl{\omega^2}^{n_2})$:
\ben
&& \rank \bE^\vee_0 =g, \\
&& \rank \bE^\vee_1 =g-1 + \frac{n_1+2n_2}{3}, \\
&& \rank \bE^\vee_2 =g-1 + \frac{2n_1+n_2}{3}.
\een
In particular,
on $\Mbar_{g,0}(\cB\bZ_3;\cl\omega^{3m})$,
$R^1\pi_*f^*V_1$ is a vector bundle of rank $g-1+m$.

By (\ref{eqn:F0}) we get
\ben
& &\sum_{n \geq 3}
\frac{1}{n!}\int_{\Mbar_{0,n}(\cB \bZ_3)} \ch_{1,1} \prod_{i=1}^n
\sum_{k=1}^{\infty} \ev_i^*(\sum_{\cl\gamma} t^{\cl\gamma}_ke_{\cl\gamma}) \bpsi_i^k  \\
& = & \frac{1}{36} (t^1_0)^2(t^2_0)^2
+  \frac{1}{216}(t^1_0)^4t^2_1+\frac{1}{216}(t^2_0)^4t^1_1
+\frac{1}{18}(t^1_0)^2(t^2_0)^2t^0_1+\frac{1}{18}t^0_0(t^1_0)^2t^2_0t^2_1 \\
&& +\frac{1}{18}t^0_0t^1_0(t^2_0)^2t^1_1+\frac{1}{27}t^1_0(t^2_0)^3t^2_1
+\frac{1}{27}(t^1_0)^3t^2_0t^1_1 + \cdots
\een
and
\ben
&&\sum_{n > 0}
\frac{1}{n!}\int_{\Mbar_{1,n}(\cB \bZ_3)} \ch_{1,1} \prod_{i=1}^n \sum_{k=1}^{\infty} \ev_i^*(t_k) \bpsi_i^k  \\
& = & \frac{1}{72} t^0_0
+ \frac{1}{24} t^1_0t^2_1 + \frac{1}{72} t^0_0t^0_1 + \frac{1}{24} t^2_0t^1_1
+ \frac{5}{144} (t^1_0)^2 t^1_2 + \frac{5}{72} t^1_0 t^2_0 t^0_2 \\
& + & \frac{1}{12} t^0_0 t^1_1 t^2_1
+  \frac{13}{288}t^0_0 t^1_0 t^2_2
+ \frac{1}{12} t^1_0 t^0_1 t^2_1
+ \frac{1}{12} t^2_0 t^0_1 t^1_1
+ \frac{1}{72} t^0_0 (t^0_1)^2 \\
& + & \frac{1}{144} (t^0_0)^2  t^0_2 + \frac{7}{192} (t^2_0)^2  t^2_2
+ \frac{1}{18} t^1_0 (t^1_1)^2  + \frac{1}{18} t^2_0 (t^2_1)^2
         + \frac{1}{24} t^0_0 t^2_0 t^1_2 + \cdots
\een
These are obtained using Maple.

Consider the action of $\bZ_3$ on $\bC^3$ where $\omega$ by multiplication by $\xi_3$.
The twisted sector of the orbifold $[\bC^3/\bZ_3]$ has two components,
corresponding to $\omega$ and $\omega^2$,
each consists of a single point, the origin.
For the first,
because $\omega$ acts on the normal space as
$\diag(e^{2\pi i /3}, e^{2\pi i/3}, e^{2\pi i/3})$,
hence the fermionic shift is
$$\frac{1}{3} + \frac{1}{3} + \frac{1}{3} = 1;$$
for the second,
because $\omega^2$ acts on the normal space as
$\diag(e^{4\pi i /3}, e^{4\pi i/3}, e^{4\pi i/3})$,
hence the fermionic shift is
$$\frac{2}{3} + \frac{2}{3} + \frac{2}{3} = 2.$$
Because the twisted degrees \cite{Bry-Gra} arise from twisted sectors of shifted fermionic
degree $1$,
to verify the predictions of Aganagic-Bouchard-Klemm \cite{Aga-Bou-Kle},
one needs to compute
\ben
\sum_m \frac{\sigma^{3m}}{3m!} \int_{\Mbar_{g, 3m}(\cB \bZ_3; \cl\omega^{3m})}
e(\bE^\vee_1 \oplus \bE^\vee_1 \oplus \bE^\vee_1).
\een
We have checked some cases using a partly developed Maple program
which automates the procedures similar to the examples in the next section.

\section{Example: The Case of $\cB \bZ_5$}

\subsection{Orbifold cohomology of $\cB Z_5$}

Let $\omega$ be a generator of $\bZ_5$.
In the class basis,
the product on the orbifold cohomology of $\cB \bZ_5$ is given by:
\be
e_{\cl{\omega^a}}e_{\cl{\omega^b}} =
 e_{\cl{\omega^{a+b}}},
\ee
and the metric is given by:
\be
\eta(e_{\cl{\omega^a}}, e_{\cl{\omega^b}}) = \begin{cases}
\frac{1}{5}, & a+ b \cong 0 \pmod{5}, \\
0, & \text{otherwise}.\end{cases}
\ee
I.e.,
\be
(\eta_{\alpha\beta})
= \begin{pmatrix}
\frac{1}{5} & 0 & 0 & 0 & 0 \\ 0 & 0 & 0 & 0 & \frac{1}{5} \\
0 & 0 & 0 & \frac{1}{5} & 0 \\ 0 & 0 & \frac{1}{5} & 0 & 0 \\
0 & \frac{1}{5} & 0 & 0 & 0
\end{pmatrix}
\ee
and so
\be
(\eta^{\alpha\beta})
= \begin{pmatrix}
5 & 0 & 0 & 0 & 0 \\ 0 & 0 & 0 & 0 & 5 \\ 0 & 0 & 0 & 5 & 0 \\
0 & 0 & 5 & 0 & 0 \\ 0 & 5 & 0 & 0 & 0
\end{pmatrix}
\ee

Let $\{ V_\alpha\}_{\alpha=0}^4$ be the set of irreducible
representations of $\bZ_5$,
where $\omega$ acts on $V_{\alpha}$ as multiplication by $\xi_5^\alpha$.
The representation basis is given by:
\[
f_\alpha := \frac{1}{5}\sum_{a=0}^4 \xi_5^{-a\alpha}e_{\cl{\omega^a}}.
\]
One has
\be
e_{\cl{\omega^a}} =  \sum_{\alpha=0}^4 \xi_5^{a\alpha} f_{\alpha}.
\ee
Let $\{t^0, t^1,\dots, t^4\}$ be linear coordinates in
$\{e_{\cl1}, e_{\cl\omega}, \dots, e_{\cl{\omega^4}}\}$
and let $\{u^0, u^1, \dots, u^4\}$ be linear coordinates in $\{f_0, f_1,\dots, f_4\}$.
\begin{align}
t^a & = \frac{1}{5}\sum_{\alpha=0}^4 \xi_5^{-a\alpha} u^\alpha, &
u^{\alpha} & = \sum_{a=0}^4 \xi_5^{a\alpha} t^a.
\end{align}

\subsection{Some Hurwitz-Hodge integrals}
Because of the monodromy condition (\ref{eqn:Monodromy}),
$\Mbar_{g,n}(\cB \bZ_5; \cl1^{n_0}, \cl\omega^{n_1}, \dots, \cl{\omega^4}^{n_4})$ is nonempty only if
$$n_1 + 2  n_2  + 3n_3 + 4n_4 \equiv 0 \pmod{5}.$$
Let $\bE^\vee_a = R^1\pi_*f^*V_a$.
By the rank formula (\ref{eqn:Rank}),
we have on the component $\Mbar_{g,n}(\cB \bZ_5; \cl1^{n_0}, \cl\omega^{n_1}, \dots, \cl{\omega^2}^{n_4})$:
\ben
&& \rank \bE^\vee_0 =g, \\
&& \rank \bE^\vee_1 =g-1 + \frac{n_1+2n_2+3n_3+4n_4}{5}, \\
&& \rank \bE^\vee_2 =g-1 + \frac{2n_1+4n_2+n_3+3n_4}{5}, \\
&& \rank \bE^\vee_3 =g-1 + \frac{3n_1+n_2+4n_3+2n_4}{5}, \\
&& \rank \bE^\vee_4 =g-1 + \frac{4n_1+3n_2+2n_3+n_4}{5}.
\een
In particular,
on $\Mbar_{g,0}(\cB\bZ_3;\cl\omega^{3m})$,
$R^1\pi_*f^*V_1$ is a vector bundle of rank $g-1+m$.

Consider the action of $\bZ_5$ on $\bC^3$ where $\omega$ acts by the diagonal matrix
$\diag (\xi_5, \xi_5, \xi_5^3)$.
The twisted sector of the orbifold $[\bC^3/\bZ_5(1,1,3)]$ has $4$ components,
corresponding to $\omega^a$, $a=1, \dots, 4$,
each consists of a single point, the origin.
Their fermionic degrees are $1$, $1$, $2$, $2$, respectively.
Because the twisted degrees \cite{Bry-Gra} arise from twisted sectors of shifted fermionic
degree $1$,
we will focus on $\Mbar_{g, n_1+n_2}(\cB \bZ_5; \cl\omega^{n_1}, \cl{\omega^2}^{n_2})$,
where $n_1+2n_2 \equiv 0 \pmod{5}$.
This space is of dimension $3g-3 + n_1 +n_2$,
and on it we have
\ben
\rank (\bE^{\vee}_1 \oplus \bE^{\vee}_1 \oplus \bE^\vee_3) = 3g-3 + n_1 + n_2.
\een
We are interested in
\ben
\sum_{n_1,n_2} \frac{\sigma_1^{n_1}\sigma_2^{n_2}}{n_1!n_2!} \int_{\Mbar_{g, n_1+n_2}(\cB \bZ_5; \cl\omega^{n_1}, \cl{\omega^2}^{n_2})}
e(\bE^\vee_1 \oplus \bE^\vee_1 \oplus \bE^\vee_3).
\een

When $n_1 + 2 n_2 = 5$,
the relevant integrals are
\bea
&& \int_{\Mbar_{0, 3}(\cB \bZ_5; \cl\omega, \cl{\omega^2}^2)} 1 = \frac{1}{5}, \label{eqn:Int1} \\
&& \int_{\Mbar_{0, 4}(\cB \bZ_5; \cl\omega^3, \cl{\omega^2})} c_1(\bE^\vee_3) = - \frac{1}{25}, \label{eqn:Int2}\\
&& \int_{\Mbar_{0, 5}(\cB \bZ_5; \cl\omega^5)} c_2(\bE^\vee_3) = \frac{1}{25}. \label{eqn:Int3}
\eea

We present some the computations to illustrate the details of the recursions.
To compute (\ref{eqn:Int2}),
note
\ben
&& \left(\frac{A_{k+1}(E_3)z^{k}}{(k+1)!} \right)^\wedge  \\
& = & - \sum_{l=0}^{\infty} \big(\frac{B_{k+1}}{(k+1)!} q_l^0\pd_{0, k+l}
+ \frac{B_{k+1}(3/5)}{(k+1)!} q_l^1\pd_{1, k+l}
+ \frac{B_{k+1}(1/5)}{(k+1)!} q_l^2\pd_{2, k+l} \\
&& + \frac{B_{k+1}(4/5)}{(k+1)!} q_l^3\pd_{3, k+l}
+ \frac{B_{k+1}(2/5)}{(k+1)!} q_l^4\pd_{4\mathrm{}, k+l}\big) \\
& + & \frac{1}{2} \hbar
\big(5\frac{B_{k+1}}{(k+1)!} \sum_{l=0}^{k-1} (-1)^l \pd_{0,l}\pd_{0,k-1-l} \\
& + & 5 \frac{B_{k+1}(2/5)}{(k+1)!} \sum_{l=0}^{k-1} (-1)^l \pd_{1,l} \pd_{4,k-1-l}
+ 5 \frac{B_{k+1}(4/5)}{(k+1)!} \sum_{l=0}^{k-1} (-1)^l \pd_{2,l} \pd_{3,k-1-l} \\
& + & 5 \frac{B_{k+1}(1/5)}{(k+1)!} \sum_{l=0}^{k-1} (-1)^l \pd_{3,l} \pd_{2,k-1-l}
+ 5 \frac{B_{k+1}(3/5)}{(k+1)!} \sum_{l=0}^{k-1} (-1)^l \pd_{4,l} \pd_{1,k-1-l}\big).
\een
Therefore,
\ben
&& \int_{\Mbar_{0,4}(\cB \bZ_5; \cl{\omega}^3, \cl{\omega^2})} \ch_{1, 3} \\
& = & \frac{B_2}{2!} \int_{\Mbar_{0,5}(\cB \bZ_5; \cl{\omega}^3, \cl{\omega^2}, \cl1)} \bpsi_5^2 \\
& - & 3 \frac{B_2(3/5)}{2!} \int_{\Mbar_{0,4}(\cB \bZ_5; \cl{\omega}^3, \cl{\omega^2})} \bpsi_1
- \frac{B_2(1/5)}{2!} \int_{\Mbar_{0,4}(\cB \bZ_5; \cl{\omega}^3, \cl{\omega^2})} \bpsi_4 \\
& + & 3 \cdot \frac{B_2(4/5)}{2!}
\int_{\Mbar_{0,3}(\cB \bZ_5; \cl{\omega}^2, \cl{\omega^3})} 1 \cdot 5 \cdot
\int_{\Mbar_{0,3}(\cB \bZ_5; \cl{\omega^2}, \cl{\omega^2}, \cl{\omega})} 1 \\
& = & \frac{B_2}{2!} \frac{1}{5}
- 3 \frac{B_2(3/5)}{2!} \frac{1}{5}
- \frac{B_2(1/5)}{2!} \frac{1}{5}
+ 3 \cdot \frac{B_2(4/5)}{2!} \cdot \frac{1}{5} \cdot 5 \cdot \frac{1}{5} \\
& = & \frac{1}{25}.
\een
In the same fashion we have
\ben
&& \int_{\Mbar_{0, 5}(\cB \bZ_5; \cl\omega^5)} \ch_{2,3} \\
& = & \frac{B_3}{3!} \int_{\Mbar_{0,6}(\cB \bZ_5; \cl{\omega}^5, \cl1)} \bpsi_6^3
- 5 \frac{B_3(3/5)}{3!} \int_{\Mbar_{0,5}(\cB \bZ_5; \cl{\omega}^5)} \bpsi_1^2 \\
& + & \binom{5}{3} \cdot \frac{B_3(1/5)}{3!}
\int_{\Mbar_{0,4}(\cB \bZ_5; \cl{\omega}^3, \cl{\omega^2})} \bpsi_+ \cdot 5 \cdot
\int_{\Mbar_{0,3}(\cB \bZ_5; \cl{\omega^3}, \cl{\omega}^2)} 1 \\
& = & \frac{B_3}{3!} \frac{1}{5}
- 5 \frac{B_3(3/5)}{3!} \frac{1}{5}
+ \binom{5}{3} \cdot \frac{B_3(1/5)}{3!} \cdot \frac{1}{5} \cdot 5 \cdot \frac{1}{5} \\
& = & \frac{1}{50}.
\een
We also have
\ben
&& \int_{\Mbar_{0,5}(\cB \bZ_5; \cl{\omega}^5)} \ch_{1,3}^2 \\
& = & \frac{B_2}{2!} \int_{\Mbar_{0,6}(\cB \bZ_5; \cl{\omega}^5, \cl1)} \ch_1(\bE^\vee_3)\bpsi_6^2 \\
& - & 5 \frac{B_2(3/5)}{2!} \int_{\Mbar_{0,5}(\cB \bZ_5; \cl{\omega}^5)} \ch_1(\bE^\vee_3) \bpsi_1 \\
& + & \binom{5}{3} \cdot \frac{B_2(1/5)}{2!}
\int_{\Mbar_{0,4}(\cB \bZ_5; \cl{\omega}^3, \cl{\omega^2})} \ch_1(\bE^\vee_3) \cdot 5 \cdot
\int_{\Mbar_{0,3}(\cB \bZ_5; \cl{\omega^3}, \cl{\omega}, \cl{\omega})} 1 \\
& = & \frac{B_2}{2!} \frac{1}{5}
- 5 \frac{B_2(3/5)}{2!} \frac{3}{25}
+ \binom{5}{3} \cdot \frac{B_2(1/5)}{2!} \cdot \frac{1}{25} \cdot 5 \cdot \frac{1}{5} \\
& = & \frac{1}{25}.
\een
In the second to last equality we have used the results of the following computations.
\ben
&& \int_{\Mbar_{0,5}(\cB \bZ_5; \cl{\omega}^5)} \ch_{1,3} \bpsi_1 \\
& = & \frac{B_2}{2!} \int_{\Mbar_{0,6}(\cB \bZ_5; \cl{\omega}^5, \cl1)} \bpsi_1\bpsi_6^2 \\
& - & \frac{B_2(3/5)}{2!} \int_{\Mbar_{0,5}(\cB \bZ_5; \cl{\omega}^5)} \bpsi_1^2
- 4 \frac{B_2(3/5)}{2!} \int_{\Mbar_{0,5}(\cB \bZ_5; \cl{\omega}^5)} \bpsi_1\bpsi_2 \\
& + & \binom{4}{2} \cdot \frac{B_2(1/5)}{2!}
\int_{\Mbar_{0,4}(\cB \bZ_5; \cl{\omega}^3, \cl{\omega^2})} \bpsi_1 \cdot 5 \cdot
\int_{\Mbar_{0,3}(\cB \bZ_5; \cl{\omega^3}, \cl{\omega}, \cl{\omega})} 1 \\
& = & \frac{B_2}{2!} \frac{3}{5}
- \frac{B_2(3/5)}{2!} \frac{1}{5}
- 4 \frac{B_2(3/5)}{2!} \frac{2}{5}
+ \binom{4}{2} \cdot \frac{B_2(1/5)}{2!} \cdot \frac{1}{5} \cdot 5 \cdot \frac{1}{5} \\
& = & \frac{3}{25}.
\een

\ben
&& \int_{\Mbar_{0,6}(\cB \bZ_5; \cl{\omega}^5)} \ch_{1,3} \bpsi_6^2 \\
& = & \frac{B_2}{2!} \int_{\Mbar_{0,7}(\cB \bZ_5; \cl{\omega}^5, \cl1^2)} \bpsi_6^2\bpsi_7^2 \\
& - & 5 \frac{B_2(3/5)}{2!} \int_{\Mbar_{0,6}(\cB \bZ_5; \cl{\omega}^5, \cl1)} \bpsi_1 \bpsi_6^2
- \frac{B_2}{2!} \int_{\Mbar_{0,6}(\cB \bZ_5; \cl{\omega}^5, \cl1)} \bpsi_6^3 \\
& + & \binom{5}{3} \cdot \frac{B_2(1/5)}{2!}
\int_{\Mbar_{0,4}(\cB \bZ_5; \cl{\omega}^3, \cl1, \cl{\omega^2})} \bpsi_6^2 \cdot 5 \cdot
\int_{\Mbar_{0,3}(\cB \bZ_5; \cl{\omega^3}, \cl{\omega}, \cl{\omega})} 1 \\
& = & \frac{B_2}{2!} \frac{6}{5}
- 5 \frac{B_2(3/5)}{2!} \frac{3}{5}
- \frac{B_2}{2!} \frac{1}{5}
+ \binom{5}{3} \cdot \frac{B_2(1/5)}{2!} \cdot \frac{1}{5} \cdot 5 \cdot \frac{1}{5} \\
& = & \frac{1}{5}.
\een
Finally,
\ben
&& \int_{\Mbar_{0,5}(\cB \bZ_5; \cl{\omega}^5)} e(\bE^\vee_3)
= \int_{\Mbar_{0,5}(\cB \bZ_5; \cl{\omega}^5)} (\half  \ch_1(\bE^\vee_3)^2 -  \ch_2(\bE^\vee_3)) \\
& = & \int_{\Mbar_{0,5}(\cB \bZ_5; \cl{\omega}^5)} (\half  \ch_{1,3} +  \ch_{2,3})
= \frac{1}{25}.
\een

{\em Acknowledgements}.
This research is partly supported by two NSFC grants and a 973 project grant.
The author also thanks Professor Youjin Zhang for providing with him the Maple V Learning Guide
and the Maple V Programming Guide.

\end{document}